\documentclass[12pt]{amsart}
\usepackage[centertags]{amsmath}
\usepackage{amsfonts,amssymb}
\usepackage{graphicx}
\usepackage{enumerate}
\usepackage[latin1]{inputenc}
\usepackage{float}
\usepackage{colonequals}

  \newtheorem{theorem}{Theorem}

  \newtheorem{lemma}{Lemma}%
  \theoremstyle{remark}

 \newcommand{\li}{{\rm li}}

\begin{document}

\title[An inequality related to the sieve of Eratosthenes]{An inequality related to\\ the sieve of Eratosthenes}
\date{\today}

\author{Kai (Steve) Fan}
\address{Mathematics Department, Dartmouth College, Hanover, NH 03755}
\email{Steve.Fan.GR@dartmouth.edu}
\author{Carl Pomerance}
\address{Mathematics Department, Dartmouth College, Hanover, NH 03755}
\email{carlp@math.dartmouth.edu}

\begin{abstract}
Let $\Phi(x,y)$ denote the number of integers $n\in[1,x]$ free of prime
factors $\le y$.  We show that but for a few small cases, $\Phi(x,y)<.6x/\log y$
when $y\le\sqrt{x}$.
\end{abstract}


\maketitle

\section{Introduction}
The sieve of Eratosthenes removes the multiples of the primes $p\le y$
from the set of positive integers $n\le x$.  Let $\Phi(x,y)$ denote
the number of integers remaining.  Answering a question of Ford, the first-named author 
\cite{F1} recently proved
the following theorem.

\medskip
\noindent{\bf Theorem A}.  {\it When $2\le y\le x$, $\Phi(x,y)< x/\log y$.}

\medskip
If $y>\sqrt{x}$, then $\Phi(x,y)=\pi(x)-\pi(y)+1$ (where $\pi(t)$ is the number
of primes in $[1,t]$), and so by the prime number theorem, Theorem A is
essentially best possible when $x^{1-\epsilon}<y<\epsilon x$.  When $y\le\sqrt{x}$, there is 
a long history in estimating $\Phi(x,y)$, and in particular, we have the following theorem, essentially due to Buchstab (see
\cite[Theorem III.6.4]{T}).

\medskip
\noindent{\bf Theorem B}. {\it For $\omega(u)$ the Buchstab function and $u=\log x/\log y\ge2$ and $y\ge2$,}
\[
\Phi(x,y)=\frac{x}{\log y}\left(\omega(u)+O\left(\frac{1}{\log y}\right)\right).
\]
 
 \medskip
The Buchstab function $\omega(u)$
is defined as the unique continuous function on $[1,\infty)$ such that 
\[
u\omega(u)=1 \hbox{ on }[1,2],\quad (u\omega(u))'=\omega(u-1)\hbox{ on }(2,\infty).
\]
Below is a graph of $\omega(u)$ for $u\in[1,8]$ generated by Mathematica. 
\begin{figure}[H]
	\includegraphics{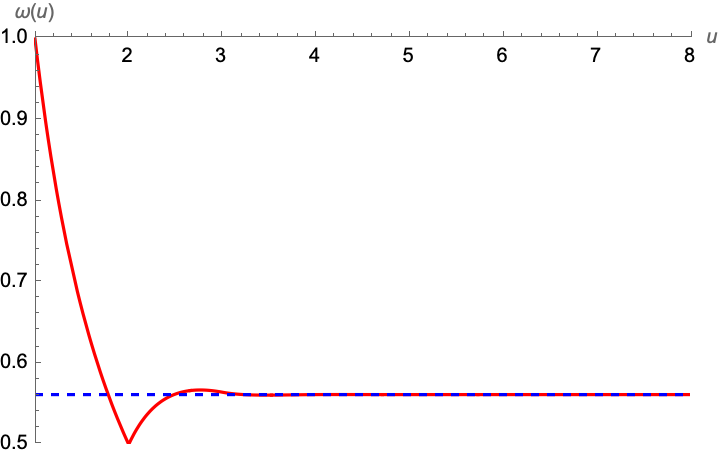}
\end{figure}
It is known that $\lim_{u\to\infty}\omega(u)=e^{-\gamma}=0.561459483566885\dots$
and that $\omega(u)$ oscillates above and below its limiting value infinitely often.
The minimum value of $\omega(u)$ on $[2,\infty)$ is $1/2$ at $u=2$ and the 
maximum value $M_0$ is $0.567143290409783\dots$,
occurring at $u=2.76322283417162\dots$.  In particular, it follows from Theorem B
that if $c >M_0$ and $y\le\sqrt{x}$ with $y$ sufficiently large depending on the choice
of $c$, that $\Phi(x,y)<cx/\log y$.  In addition, using an inclusion--exclusion argument
plus the fact that the Mertens product $\prod_{p\le y}(1-1/p)<M_0/\log y$ for all $y\ge2$,
the inequality $\Phi(x,y)<cx/\log y$ can be extended to all $2\le y\le\sqrt{x}$, but now
with $x$ sufficiently large depending on $c$.

In light of Theorem A and given that $\Phi(x,y)$ is a fundamental (and ancient) function,
 it seems interesting
to try and make these consequences of Theorem B numerically explicit.  We prove the following theorem.
\begin{theorem}
\label{th:ineq}
For $3\le y\le\sqrt{x}$, we have $\Phi(x,y)<.6x/\log y$.  The same inequality holds
when $2\le y\le\sqrt{x}$ and $x\ge10$.
\end{theorem}

To prove this we use some numerically explicit estimates of primes due to
Rosser--Schoenfeld,  B\"uthe, and others.
 In addition we use a numerically
explicit version of the upper bound in Selberg's sieve.  

Theorem B itself is also appealing. It provides a simple asymptotic formula for $\Phi(x,y)$ as $y\to\infty$ which is applicable in a wide range. Writing 
\[\Phi(x,y)=\frac{x}{\log y}\left(\omega(u)+\frac{\Delta(x,y)}{\log y}\right),\]
one may attempt to establish numerically explicit lower and upper bounds for $\Delta(x,y)$ in the range $y\le\sqrt{x}$ for suitably large $y\ge y_0$, where $y_0\ge2$ is some numerically computable constant. More precisely, de Bruijn \cite{Br} essentially showed that for any given $\epsilon>0$, one has
\[\Phi(x,y)=\mu_y(u)e^{\gamma}x\log y\prod_{p\leq y}\left(1-\frac{1}{p}\right)+O(x\exp(-(\log y)^{3/5-\epsilon}))\]
for all $x\ge y\ge2$, where 
\[\mu_y(u)\colonequals\int_{0}^{u-1}\omega(u-v)y^{-v}\,dv.\]
Recently, the first-named author \cite{F2} proved numerically explicit versions of this result applicable for $y$ in wide ranges.

\section{A prime lemma}

Let $\pi(x)$ denote, as usual, the number of primes $p\le x$.
Let 
\[
\li(x)=\int_0^x\frac{dt}{\log t},
\]
where the principal value is taken for the singularity at $t=1$.
There is a long history in trying to find the first point when $\pi(x)\ge\li(x)$,
which we now know is beyond $10^{19}$.  We prove a lemma
based on this research.
\begin{lemma}
\label{lem:pnt}
Let $\beta_0=2.3\times10^{-8}$.
For $x\ge 2$, we have $\pi(x)<(1+\beta_0)\li(x)$.
\end{lemma}
\begin{proof}
The result is true for $x\le 10$, so assume $x\ge10$.
Consider the Chebyshev function 
\[
\theta(x)=\sum_{p\le x}\log p.
\]
We use \cite[Prop. 2.1]{LP}, which depends strongly on extensive calculations
of B\"uthe \cite{B1,B2} and Platt \cite{P}.
This result asserts in part that $\theta(x)\le x-.05\sqrt{x}$ for $1427\le x\le10^{19}$ and for
larger $x$, $\theta(x)<(1+\beta_0)x$.  One easily checks that $\theta(x)<x$ for $x<1427$,
so we have
\[
\theta(x)<(1+\beta_0)x,\quad x>0.
\]
By partial summation, we have
\begin{align*}
\pi(x)&=\frac{\theta(x)}{\log x}+\int_2^x\frac{\theta(t)}{t(\log t)^2}\,dt\\
&<\frac{(1+\beta_0)x}{\log x}+\int_2^{10}\frac{\theta(t)}{t(\log t)^2}\,dt+(1+\beta_0)\int_{10}^x\frac{dt}{(\log t)^2}.
\end{align*}
Since $\int dt/(\log t)^2 =-t/\log t+\li(t)$, we have
\begin{align}
\label{eq:pntub}
\pi(x)& <(1+\beta_0)\li(x) +\int_2^{10}\frac{\theta(t)}{t(\log t)^2}\,dt+(1+\beta_0)(10/\log 10-\li(10))\notag\\
&<(1+\beta_0)\li(x) -.144.
\end{align}
This gives the lemma.
\end{proof}
After checking for $x\le10$, we remark that an immediate corollary of \eqref{eq:pntub} is the inequality
\begin{equation}
\label{eq:pntub2}
\pi(x)-k<(1+\beta_0)(\li(x)-k),\quad 2\le k\le \pi(x),~ k\le 10^7.
\end{equation}

\section{Inclusion--exclusion}
\label{sec:incexc}
For small values of $y\ge2$ we can do a complete inclusion--exclusion to compute
$\Phi(x,y)$.  Let $P(y)$ denote the product of the primes $p\le y$.
We have
\begin{equation}
\label{eq:inex}
\Phi(x,y)=\sum_{d\mid P(y)}\mu(d)\left\lfloor\frac xd\right\rfloor.
\end{equation}
As a consequence, we have
\begin{equation}
\label{eq:elub}
\Phi(x,y)\le\sum_{d\mid P(y)}\mu(d)\frac xd+\sum_{\substack{d\mid P(y)\\\mu(d)=1}}1
=x\prod_{p\le y}\left(1-\frac1p\right)+2^{\pi(y)-1}.
\end{equation}

We illustrate how this elementary inequality can be used in the case when $\pi(y)=5$,
that is, $11\le y<13$.  Then the product in \eqref{eq:elub} is $16/77<.207793$.
The remainder term in \eqref{eq:elub} is 16.  And we have
\[
\Phi(x,y)<.207793 x+16<.6x/\log 13
\]
when $x\ge613$.  There remains the problem of dealing with smaller values of $x$,
which we address momentarily.  We apply this method for $y<71$.

\begin{footnotesize}
\begin{table}[ht]
\caption{Small $y$.}
\label{Ta:smally}
\begin{tabular}{|l|l|l|} \hline
$y$ interval & $x$ bound & max \\ \hline
~$[2,3)$ & 22 & .61035\\
~$[3,5)$ &  51 & .57940\\
~$[5,7)$ &  96 & .55598\\
~$[7,11)$ &  370 &.56634\\
~$[11,13)$ & 613 & .55424\\
~$[13,17)$ &  1603 & .56085\\
~$[17,19)$ & 2753 & .54854 \\
~$[19,23)$ &  6296 & .55124\\
~$[23,29)$ &  17539 & .55806\\
~$[29,31)$ & 30519  & .55253\\
~$[31,37)$ & 76932 & .55707\\
~$[37,41)$ & $1.6\times10^5$ & .55955\\
~$[41,43)$ & $2.9\times10^5$ & .55648\\
~$[43,47)$ & $5.9\times10^5$ & .55369\\
~$[47,53)$ & $1.4\times10^6$ & .55972\\
~$[53,59)$ & $3.0\times 10^6$ & .55650\\
~$[59,61)$ &  $5.4\times10^6$ & .55743\\
~$[61,67)$ &  $1.2\times10^7$ & .55685\\
~$[67,71)$ & $2.4\times10^7$ & .55641\\
  \hline
   \end{tabular}
\end{table}
\end{footnotesize}

Pertaining to Table \ref{Ta:smally}, for $x$ beyond the ``$x$ bound" and $y$ in the given interval, we have
$\Phi(x,y) < .6x/\log y$.
The column ``max" in Table \ref{Ta:smally} is the supremum of $\Phi(x,y)/(x/\log y)$
for $y$ in the given interval and $x\ge y^2$ with $x$ below the $x$ bound.  
The max statistic was computed by creating a table
of the integers up to the $x$ bound with a prime factor $\le y$, taking the complement
of this set in the set of all integers up to the $x$ bound, and then computing $(j\log p)/n$
where $n$ is the $j$th member of the set and $p$ is the upper bound of the $y$ interval.
The max of these numbers is recorded as the max statistic.

As one can see, for $y\ge 3$  the max statistic in Table \ref{Ta:smally} is
below $.6$.  However,
for the interval $[2,3)$ it is above $.6$.  One can compute
that it is $<.6$ once $x\ge10$. 

This method can be extended to larger values of $y$, but the $x$ bound becomes
prohibitively large.  With a goal of keeping the $x$ bound smaller than $3\times10^7$,
we can extend a version of inclusion-exclusion to $y<241$ as follows.

First, we ``pre-sieve" with the primes 2, 3, and 5.  For any $x\ge0$ the number of integers
$n\le x$ with $\gcd(n,30)=1$ is $(4/15)x+r$, where $|r|\le14/15$, as can be easily verified
by looking at values of $x\in[0,30]$.  We change the definition of $P(y)$ to be the product
of the primes in $(5,y]$.  Then for $y\ge5$, we have
\[
\Phi(x,y)\le\frac4{15}\sum_{d\mid P(y)}\mu(d)\frac xd+\frac{14}{15}2^{\pi(y)-3}.
\]
However, it is better to use the Bonferroni inequalities in the form
\[
\Phi(x,y)\le \frac4{15}\sum_{j\le4}\sum_{\substack{d\mid P(y)\\\nu(d)=j}}(-1)^j\frac xd
+
\sum_{i=0}^4\binom{\pi(y)-3}{i}
=xs(y)+b(y),
\]
say, where $\nu(d)$ is the number of distinct prime factors of $d$.
(We remark that the expression $b(y)$ could be replaced with $\frac{14}{15}b(y)$.)
  The inner sums in $s(y)$ can be computed easily
using Newton's identities, and we see that
\[
\Phi(x,y)\le.6x/\log y \hbox{ for }x>b(y)/(.6/\log y-s(y)).
\]
We have verified that this  $x$ bound is smaller than $30{,}000{,}000$ for $y<241$
and we have verified that $\Phi(x,y) <.6x/\log y$ for $x$ up to this bound and $y<241$.

This completes the proof of Theorem \ref{th:ineq} for $y<241$.
 
\section{When $u$ is large:  Selberg's sieve}
\label{sec:ularge}
In this section we prove Theorem \ref{th:ineq} in the case that $u=\log x/\log y\ge7.5$
and $y\ge241$.  Our principal tool is a numerically explicit form of Selberg's sieve.

Let $\mathcal A$ be a set of positive integers $a\le x$ and with $|\mathcal A|\approx X$.
  Let
$\mathcal P=\mathcal P(y)$ be a set of primes $p\le y$.  For each $p\in\mathcal P$
we have a collection of $\alpha(p)$ residue classes mod $p$, where $\alpha(p)<p$.
Let $P=P(y)$ denote the product of the members of $\mathcal P$.  Let $g$ be the
multiplicative function defined for numbers $d\mid P$ where $g(p)=\alpha(p)/p$
when $p\in\mathcal P$.  We let
\[
V:=\prod_{p\in\mathcal P}(1-g(p))=\prod_{p\in\mathcal P}\left(1-\frac{\alpha(p)}p\right).
\]
We define $r_d(\mathcal A)$ via the equation
\[
\sum_{\substack{a\in\mathcal A\\d\mid a}}1=g(d)X+r_d(\mathcal A).
\]
The thought is that $r_d(\mathcal A)$ should be small.
We are interested in 
$S(\mathcal A,\mathcal P)$, the number of those $a\in\mathcal A$ such that $a$ is coprime to $P$.

We will use Selberg's sieve as given in \cite[Theorem 7.1]{FI}.  This involves an auxiliary
parameter $D<X$ which can be freely chosen.  Let $h$ be the multiplicative function supported
on divisors of $P$ such that $h(p)=g(p)/(1-g(p))$.  In particular if each $\alpha(p)=1$, then
each $g(p)=1/p$ and $h(p)=1/(p-1)$, so $h(d)=1/\varphi(d)$ for $d\mid P$, where
$\varphi$ is Euler's function.  Henceforth we will make this assumption (that each
$\alpha(p)=1$).  Let
\[
J=J_D=\sum_{\substack{d\mid P\\d<\sqrt{D}}}h(d),\quad
R=R_D=\sum_{\substack{d\mid P\\d<D}}\tau_3(d)|r_d(\mathcal A)|,
\]
where $\tau_3(d)$ is the number of ordered factorizations $d=abc$, where $a,b,c$ 
are positive integers.  Selberg's sieve gives in this situation that
\begin{equation}
\label{eq:selberg}
S(\mathcal A, \mathcal P)\le X/J+R.
\end{equation}
Note that if $D\ge P^2$, then
\[
J=\sum_{d\mid P}h(d)=\prod_{p\in\mathcal P}(1+h(p))=\prod_{p\in\mathcal P}(1-g(p))^{-1} = V^{-1},
\]
so that $X/J=XV$.  This is terrific, but if $D$ is so large, the remainder term $R$ in \eqref{eq:selberg}
is also large, making the estimate useless.  So, the trick is to choose $D$ judiciously
so that $R$ is under control with $J$ being near to $V^{-1}$.

Consider the case when each $|r_d(\mathcal A)|\le r$, for a constant $r$.
In this situation the following lemma is useful.
\begin{lemma}
\label{lem:R}
For $y\ge241$, we have
\[
R\le r\sum_{\substack{d<D\\d\mid P(y)}}\tau_3(d)\le rD(\log y)^2\prod_{\substack{p\le y\\p\notin\mathcal P}}
\left(1+\frac2p\right)^{-1}.
\]
\end{lemma}
\begin{proof}
Let $\tau(d)=\tau_2(d)$ denote the number of positive divisors of $d$.  Note that
\[
\sum_{d\mid P(y)}\frac{\tau(d)}d=
\prod_{p\in\mathcal P}\left(1+\frac2p\right)
=\prod_{p\le y}\left(1+\frac2p\right)\prod_{\substack{p\le y\\p\notin\mathcal P}}\left(1+\frac2p\right)^{-1}.
\]
One can show that for $y\ge241$ the first product on the right  is smaller than $.95(\log y)^2$,
but we will only use the ``cleaner" bound $(\log y)^2$ (which holds when $y\ge53$).  Thus,
\begin{align*}
\sum_{\substack{d<D\\d\mid P(y)}}\tau_3(d)&=\sum_{\substack{d<D\\d\mid P(y)}}\sum_{j\mid d}\tau(j)
\le \sum_{\substack{j<D\\j\mid P(y)}}\tau(j)\sum_{\substack{d<D/j\\d\mid P(y)}}1\\
&<D\sum_{\substack{j<D\\j\mid P(y)}}
\frac{\tau(j)}j<D(\log y)^2\prod_{\substack{p\le y\\p\notin\mathcal P}}\left(1+\frac2p\right)^{-1}.
\end{align*}
This completes the proof.
\end{proof}
To get a lower bound for $J$ in \eqref{eq:selberg} we proceed as in \cite[Section 7.4]{FI}.
Recall that we are assuming  each $\alpha(p)=1$ and
so $h(d)=1/\varphi(d)$ for $d\mid P$.

Let
\[
I = \sum_{\substack{d\ge\sqrt{D}\\d\mid P}}\frac1{\varphi(d)},
\]
so that $I+J=V^{-1}$.  Hence
\begin{equation}
\label{eq:J}
J=V^{-1}-I=V^{-1}(1-IV),
\end{equation}
so we want an upper bound for $IV$.  
Let $\varepsilon$ be arbitrary with $\varepsilon>0$.  We have
\[
I<D^{-\varepsilon}\sum_{d\mid P}\frac{d^{2\varepsilon}}{\varphi(d)}
=D^{-\varepsilon}\prod_{p\le y}\left(1+\frac{p^{2\varepsilon}}{p-1}\right),
\]
and so, assuming each $\alpha(p)=1$,
\begin{equation}
\label{eq:rankin}
IV<D^{-\varepsilon}\prod_{p\in\mathcal P}\left(1+\frac{p^{2\varepsilon}-1}{p}\right)=:f(D,\mathcal P,\varepsilon).
\end{equation}
In particular, if $y\ge241$ and each $r_d(\mathcal A)\le r$, then
\begin{equation}
\label{eq:epsest}
S(\mathcal A,\mathcal P)\le XV\big(1-f(D,\mathcal P,\varepsilon)\big)^{-1}
+rD(\log y)^2\prod_{\substack{p\le y\\p\notin\mathcal P}}
\left(1+\frac2p\right)^{-1}.
\end{equation}
We shall choose $D$ so that the remainder term is small in comparison to $XV$,
and once $D$ is chosen, we shall choose $\varepsilon$ so as to minimize $f(D,\mathcal P,\varepsilon)$.

\medskip
\noindent
{\bf 4.1. The case when $y\le500{,}000$ and $u\ge7.5$.}

\medskip
\noindent
We wish to apply \eqref{eq:epsest} to estimate $\Phi(x,y)$ when $u\ge7.5$, that is, when
$x\ge y^{7.5}$.  We have a few choices for $\mathcal A$ and $\mathcal P$.  The most natural
choice is that $\mathcal A$ is the set of all integers $\le x$, $X=x$,
and $\mathcal P$ is the set of all primes $\le y$.  In this case, each $r_d(\mathcal A)\le 1$,
so that we can take $r=1$ in \eqref{eq:epsest} (since $r_d(\mathcal A)\ge 0$ in this case).
Instead we choose (as in the last section) $\mathcal A$ as the set of all integers $\le x$ that are coprime to 30
and we choose $\mathcal P$ as the set of primes $p$ with $7\le p\le y$.  Then $X=4x/15$
and one can check that each $|r_d(\mathcal A)|\le 14/15$, so we can take $r=14/15$ in
\eqref{eq:epsest}.  Also, 
\[
\prod_{\substack{p\le y\\p\notin\mathcal P}}\left(1+\frac2p\right)^{-1}=\frac3{14},
\]
when $y\ge5$.  With this choice of $\mathcal A$ and $\mathcal P$, \eqref{eq:epsest}
becomes
\begin{equation}
\label{eq:30trick}
\Phi(x,y)\le XV\left(1-D^{-\varepsilon}\prod_{7\le p\le y}\left(1+\frac{p^{2\varepsilon}-1}p\right)\right)^{-1}+\frac15D(\log y)^2,
\end{equation}
when $y\ge241$.

Our ``target" for $\Phi(x,y)$ is $.6x/\log y$.
We choose $D$ here so that our
estimate for the remainder term is 1\% of the target, namely $.006x/\log y$.
Thus, in light of Lemma \ref{lem:R}, we choose
\[
D=.03x/(\log y)^3.
\]

We have verified that for every value of $y\le500{,}000$ and $x\ge y^{7.5}$ that
the right side of \eqref{eq:30trick} is smaller than $.6x/\log y$.  Note that to verify
this, if $p,q$ are consecutive primes with $241\le p<q$, then $S(\mathcal A,\mathcal P)$
is constant for $p\le y<q$, and so it suffices to show the right side of \eqref{eq:30trick}
is smaller than $.6x/\log q$.  Further, it suffices to take $x=p^{7.5}$, since as $x$
increases beyond this point with $\mathcal P$ and $\varepsilon$ fixed, the 
expression $f(D,\mathcal P,\varepsilon)$ decreases.  For smaller values of $y$
in the range, we used Mathematica to choose the optimal choice of $\varepsilon$.
For larger values, we let $\varepsilon$ be a judicious constant over a long interval.
As an example, we chose $\varepsilon=.085$ in the top half of the range.

\medskip
\noindent
{\bf 4.2. When $y\ge 500{,}000$ and $u\ge7.5$.}
\medskip

As in the discussion above we have a few choices to make, namely for the
quantities $D$ and $\varepsilon$.  First, we choose $x= y^{7.5}$,
since the case $x\ge y^{7.5}$ follows from the proof of the case of equality.  
  We choose $D$ as before, namely $.03x/(\log y)^3$.
We also choose
\[
\varepsilon = 1/\log y.
\]
 Our goal is to prove a small
upper bound for $f(D,\mathcal P,\varepsilon)$ given in \eqref{eq:rankin}.
We have
\[
f(D,\mathcal P,\varepsilon)<D^{-\varepsilon}\exp\left(\sum_{7\le p\le y}\frac{p^{2\varepsilon}-1}p\right).
\]

We treat the two sums separately.  First, by Rosser--Schoenfeld \cite[Theorems 9, 20]{RS}, one can show that
\[
-\sum_{p\le y}\frac1p<-\log\log y-.26
\]
for all $y\ge2$, so that 
\begin{equation}
\label{eq:dusart}
-\sum_{7\le p\le y}\frac1p<-\log\log y-.26+31/30
\end{equation}
for $y\ge 7$.
For the second sum we have
\[
\sum_{7\le p\le y}p^{2\varepsilon-1}=7^{2\varepsilon-1}+(\pi(y)-4)y^{2\varepsilon-1}+
\int_{11}^y(1-2\varepsilon)(\pi(t)-4)t^{2\varepsilon-2}\,dt.
\]
At this point we use 
\eqref{eq:pntub2} , so that
\begin{align*}
\frac1{1+\beta_0}&\sum_{11\le p\le y}p^{2\varepsilon-1}
<(\li(y)-4)y^{2\varepsilon-1}
+\int_{11}^y(1-2\varepsilon)(\li(t)-4)t^{2\varepsilon-2}\,dt\\
&=(\li(y)-4)y^{2\varepsilon-1}-(\li(t)-4)t^{2\varepsilon-1}\Big|_{11}^y+\int_{11}^y\frac{t^{2\varepsilon-1}}{\log t}\,dt\\
&=(\li(11)-4)11^{2\varepsilon-1}+\li(t^{2\varepsilon})\Big|_{11}^y\\
&=(\li(11)-4)11^{2\varepsilon-1}+\li(y^{2\varepsilon})-\li(11^{2\varepsilon}),
\end{align*}
and so
\begin{align}
\label{eq:psum}
\frac1{1+\beta_0}\sum_{7\le p\le y}p^{2\varepsilon-1}
<7^{2\varepsilon-1}+(\li(11)-4)11^{2\varepsilon-1}+\li(y^{2\varepsilon})-\li(11^{2\varepsilon}).
\end{align}
There are a few things to notice, but we will not need them.  For example,
$\li(y^{2\varepsilon})=\li(e^2)$ and $\li(11^{2\varepsilon})\approx\log(11^{2\varepsilon}-1)+\gamma$.

Let $S(y)$ be the sum of the right side of \eqref{eq:dusart} and $1+\beta_0$ times
the right side of \eqref{eq:psum}.
Then
\[
f(D,\mathcal P,\varepsilon)<D^{-\varepsilon}e^{S(y)}.
\]
The expression $XV$ in \eqref{eq:30trick} is 
\[
x\prod_{p\le y}\left(1-\frac1p\right).
\]
We know from \cite{LP} that this product is $<e^{-\gamma}/\log y$ for $y\le 2\times10^9$,
and for larger values of $y$, it follows from \cite[Theorem 5.9]{D} (which proof follows
from \cite[Theorem 4.2]{D} or \cite[Corollary 11.2]{B})
that it is $<(1+2.1\times10^{-5})e^{-\gamma}/\log y$.
We have
\begin{align}
\label{eq:bigineq}
\Phi(x,y)&\le XV\big(1-f(D,\mathcal P,\varepsilon)\big)^{-1}+\frac15D(\log y)^2\\
&<(1+2.1\times10^{-5})\frac x{e^\gamma \log y}\big(1-D^{-\varepsilon}e^{S(y)}\big)^{-1}+\frac{.006x}{\log y}.\notag
\end{align}
We have verified that $(1-D^{-\varepsilon}e^{S(y)})^{-1}$ is decreasing in $y$, and that at $y=500{,}000$ it is
smaller than $1.057$.  
Thus, \eqref{eq:bigineq} implies that
\[
\Phi(x,y)<(1+2.1\times10^{-5})\frac{1.057x}{e^\gamma\log y}+\frac{.006x}{\log y}
<\frac{.5995x}{\log y}.
\]
This concludes the case of $u\ge 7.5$.

\section{Small $u$}
In this section we prove that $\Phi(x,y) < .57163x/\log y$ when $u\in[2,3)$, that
is, when $y^2\le x<y^3$.

For small values of $y$, we calculate the maximum of $\Phi(x,y)/(x/\log y)$
for $y^2\le x<y^3$ directly, as we did in Section \ref{sec:incexc} when we checked
 below the
$x$ bounds~in Table \ref{Ta:smally} and the bound $3\times10^7$.  We have done this for $241\le y\le 1100$,
and in this range we have 
\[
\Phi(x,y)<.56404\frac x{\log y},\quad y^2\le x<y^3,\quad 241\le y\le1100.
\]

Suppose now that $y>1100$ and $y^2\le x<y^3$.  We have
\begin{equation}
\label{eq:buch}
\Phi(x,y)=\pi(x)-\pi(y)+1+\sum_{y<p\le x^{1/2}}(\pi(x/p)-\pi(p)+1).
\end{equation}
Indeed, if $n$ is counted by $\Phi(x,y)$, then $n$ has at most 2 prime factors
(counted with multiplicity), so $n=1$, $n$ is a prime in $(y,x]$ or $n=pq$,
where $p,q$ are primes with $y<p\le q\le x/p$.  

Let $p_j$ denote the $j$th prime.  Note that
\[
\sum_{p\le t}\pi(p)=\sum_{j\le\pi(t)}j=\frac12\pi(t)^2+\frac12\pi(t).
\]
Thus,
\[
\sum_{y<p\le x^{1/2}}(\pi(p)-1)=\frac12\pi(x^{1/2})^2-\frac12\pi(x^{1/2})-\frac12\pi(y)^2+\frac12\pi(y),
\]
and so
\begin{equation}
\label{eq:ident}
\Phi(x,y)=\pi(x)-M(x,y)+\sum_{y<p\le x^{1/2}}\pi(x/p),
\end{equation}
where
\[
M(x,y)=\frac12\pi(x^{1/2})^2-\frac12\pi(x^{1/2})-\frac12\pi(y)^2+\frac32\pi(y)-1.
\]

We use Lemma \ref{lem:pnt} on various terms in \eqref{eq:ident}.
In particular, we have (assuming $y\ge5$)
\begin{equation}\label{eq:buthineq}
\Phi(x,y)<(1+\beta_0)\li(x)+\sum_{y<p\le x^{1/2}}(1+\beta_0)\li(x/p)
-M(x,y).
\end{equation}

Via partial summation, we have
\begin{equation}
\begin{split}
\label{eq:lixp}
\sum_{y<p\le x^{1/2}}\li(x/p)=&~x^{1/2}\li(x^{1/2})\sum_{y<p\le x^{1/2}}\frac1p\\
&~~-\int_y^{x^{1/2}}\left(\li(x/t)-\frac{x/t}{\log(x/t)}\right)\sum_{y<p\le t}\frac1p\,dt.
\end{split}
\end{equation}

For $1100\le t\le 10^4$ we have checked numerically that
\[
0<\sum_{p\le t}\frac1p-\log\log t-B<.00624,
\]
where $B=.261497\dots$ is the Meissel--Mertens constant.  Further, for
$10^4\le t\le 10^6$,
\[
0<\sum_{p\le t}\frac1p-\log\log t-B<.00161.
\]
(The lower bounds here follow as well from \cite[Theorem 20]{RS}.)
It thus follows for
$1100\le y\le10^4$ that
\begin{equation}
\label{eq:10^4}
\sum_{y<p\le x^{1/2}}\frac1p<\log\frac{\log(x^{1/2})}{\log y}+\beta_1,\quad \sum_{y<p\le t}\frac1p>\log\frac{\log t}{\log y}-\beta_1,
\end{equation}
where $\beta_1=.00624$.
Now suppose that $y\ge10^4$.  Using
\cite[Eq. (5.7)]{D} and the value 4.4916 for ``$\eta_3$" from \cite[Table 15]{B}, we have that
\[
\Big|\sum_{p\le t}\frac1p-\log\log t-B\Big|<1.9036/(\log t)^3,~~t\ge10^6.
\]
Thus, \eqref{eq:10^4} continues to hold for $y\ge10^4$ with .00624 improved to .00322.
We thus have from \eqref{eq:lixp}
\begin{equation}
\label{eq:lixp2}
\begin{split}
\sum_{y<p\le x^{1/2}}\li(x/p)<&~x^{1/2}\li(x^{1/2})\left(\log\frac{\log(x^{1/2})}{\log y}+\beta_1\right)\\
&~~-\int_y^{x^{1/2}}\left(\li(x/t)-\frac{x/t}{\log(x/t)}\right)\left(\log\frac{\log t}{\log y}-\beta_1\right)\,dt.
\end{split}
\end{equation}

Let $R(t)=(1+\beta_0)\li(t)/(t/\log t)$, so that $R(t)\to1+\beta_0$ as $t\to\infty$.  We write the first term 
on the right side of 
\eqref{eq:buthineq} as
\[
\frac x{u\log y}R(x)=\frac{R(y^u)}u\frac x{\log y},
\]
and note that the first term on the right of \eqref{eq:lixp2} is less than
\[
R(y^{u/2})\frac2u(\log(u/2)+\beta_1)\frac x{\log y}.
\]

For the expression
 $\frac12\pi(x^{1/2})^2-\frac12\pi(x^{1/2})$ in $M(x,y)$ we use the inequality $\pi(t)>t/\log t+t/(\log t)^2$
when $t\ge599$, which follows from \cite[Lemma 3.4]{BD} and a calculation
(also see \cite[Corollary 5.2]{D}).
Further, we use $\pi(y)\le R(y)y/\log y$ for the rest of $M(x,y)$.

Using these estimates and numerical integration for the integral in \eqref{eq:lixp2}
we find that 
\[
\Phi(x,y)<.57163\frac{x}{\log y},\quad y\ge1100,\quad y^2\le x<y^3.
\]

\section{Iteration}
Suppose $k$ is a positive integer and we have shown that
\begin{equation}
\label{eq:levelk}
\Phi(x,y)\le c_k\frac x{\log y}
\end{equation}
for all $y\ge241$ and $u=\log x/\log y \in[2,k)$. We can
try to find some $c_{k+1}$ not much larger than $c_k$ such that
\[
\Phi(x,y)\le c_{k+1}\frac x{\log y}
\]
for $y\ge 241$ and $u<k+1$.  We start with $c_3$, which by the results of the previous section
we can take as .57163.  In this section we attempt to find $c_k$ for $k\le 8$ such that $c_8<.6$.
It would then follow from Section \ref{sec:ularge} that $\Phi(x,y)<.6x/\log y$ for all $u\ge2$ and
$y\ge241$.

Suppose that \eqref{eq:levelk} holds and that $y$ is such that $x^{1/(k+1)}<y\le x^{1/k}$.  We have
\begin{equation}
\label{eq:iter}
\Phi(x,y)=\Phi(x,x^{1/k})+\sum_{y<p\le x^{1/k}}\Phi(x/p,p^-).
\end{equation}
Indeed the sum counts all $n\le x$ with least prime factor $p\in(y,x^{1/k}]$, and $\Phi(x,x^{1/k})$
counts all $n\le x$ with least prime factor $>x^{1/k}$.  As we have seen, it suffices to deal
with the case when $y=q_0^-$ for some prime $q_0$.

Note that if $\eqref{eq:levelk}$ holds,
then it also holds for $y=x^{1/k}$.  Indeed, if $y$ is a prime, then $\Phi(x,y)=\Phi(x,y+\epsilon)$
for all $0<\epsilon<1$, and in this  case $\Phi(x,y)\le c_kx/\log(y+\epsilon)$, by hypothesis.
Letting $\epsilon\to0$ shows we have $\Phi(x,y)\le c_kx/\log y$ as well.  If $y$ is not
prime, then for all sufficiently small $\epsilon>0$, we again have $\Phi(x,y)=\Phi(x,y+\epsilon)$
and the same proof works.

Thus, we have \eqref{eq:levelk} holding for all of the terms on the right side of \eqref{eq:iter}.
This implies that
\begin{equation}
\label{eq:ckiter}
\Phi(x,q_0^-)\le c_kx\Bigg(\frac1{\log(x^{1/k})}+\sum_{q_0\le p\le x^{1/k}}\frac1{p\log p}\Bigg).
\end{equation}
We expect that the parenthetical expression here is about the same as $1/\log q_0$, so let us try
to quantify this.  Let
\[
\epsilon_k(q_0)=\max\Bigg\{\frac{-1}{\log q_0}+\frac1{\log(x^{1/k})}+\sum_{q_0\le p\le x^{1/k}}\frac1{p\log p}:
y^k<x\le y^{k+1}\Bigg\}.
\]
Let $q_1$ be the largest prime $\le x^{1/k}$, so that
\[
\epsilon_k(q_0)=\max\Bigg\{\frac{-1}{\log q_0}+\frac1{\log q_1}+\sum_{q_0\le p\le q_1}\frac1{p\log p}:
q_0<q_1\le q_0^{1+1/k}\Bigg\}.
\]
It follows from \eqref{eq:ckiter} that
\[
\Phi(x,y)=\Phi(x,q_0^-)\le c_kx\left(\frac1{\log q_0}+\epsilon_k(q_0)\right)
=\frac{c_kx}{\log y}(1+\epsilon_k(q_0)\log q_0).
\]

Note that as $k$ grows, $\epsilon_k(q_0)$ is non-increasing since the max is over a
smaller set of primes $q_1$.  Thus, we have the inequality
\begin{equation}
\label{eq:geom}
\Phi(x,q_0^-)\le c_3(1+\epsilon_3(q_0)\log q_0)^j\frac x{\log y},\quad x^{1/3}< q_0\le x^{1/(3+j)}.
\end{equation}
Thus, we would like
\begin{equation}
\label{eq:8goal}
c_3(1+\epsilon_3(q_0)\log q_0)^5 < .6
\end{equation}

We have checked \eqref{eq:8goal} numerically for primes $q_0<1000$ and it holds
for $q_0\ge241$.

This leaves the case of primes $>1000 $.  We have the identity
\begin{align*}
\sum_{q_0\le p\le q_1}&\frac1{p\log p}\\
=&\frac{-\theta(q_0^-)}{q_0(\log q_0)^2}
+\frac{\theta(q_1)}{q_1(\log q_1)^2}+\int_{q_0}^{q_1}\theta(t)\left(\frac1{t^2(\log t)^2}+\frac2{t^2(\log t)^3}\right)\,dt,
\end{align*}
via partial summation, where $\theta$ is again Chebyshev's function.  First assume that
$q_1<10^{19}$.  Then, using \cite{B1}, \cite{B2}, we have $\theta(t)\le t$, so that
\[
\sum_{q_0\le p\le q_1}\frac1{p\log p}<\frac{q_0-\theta(q_0^-)}{q_0(\log q_0)^2}+\frac1{\log q_0}-\frac1{\log q_1}.
\]
We also have \cite{B1}, \cite{B2} that $q_0-\theta(q_0^-)<1.95\sqrt{q_0}$, so that one can
verify that
\[
\epsilon_3(q_0)<\frac{1.95}{\sqrt{q_0}(\log q_0)^2}
\]
and so \eqref{eq:8goal} holds for $q_0>1000$.  It remains to consider the cases when
$q_1>10^{19}$, which implies $q_0>10^{14}$.  Here we use $|\theta(t)-t|<3.965t/(\log t)^2$,
which is from \cite[Theorem 4.2]{D} or \cite[Corollary 11.2]{B}.  This shows that \eqref{eq:8goal} holds here as well,
completing the proof of Theorem \ref{th:ineq}.

\end{document}